\documentclass{ifacconf}

\usepackage{graphicx}      
\usepackage{natbib}        
\usepackage{amssymb,amsmath,dsfont,xcolor}

\usepackage{algorithm}
\usepackage{algpseudocode}

\newtheorem{remark}{Remark}

\begin{document}
	\begin{frontmatter}
		
		\title{Distributed Interior Point Methods for Optimization in Energy Networks} 
		
		
		\author[First]{Alexander Engelmann} 
			\author[First]{Michael Kaupmann} 
		\author[First]{and Timm Faulwasser} 
		
		\address[First]{Institute for Energy Systems, Energy Efficiency and Energy Economics, TU Dortmund University, Dortmund, Germany \\(e-mail: \{alexander.engelmann, timm.faulwasser\}@ieee.org).}
		
		\begin{abstract}                
		This note discusses an essentially decentralized interior point method, which is well suited for optimization problems arising in energy networks. 
			Advantages of the proposed method are  guaranteed and fast local convergence also for problems with non-convex constraints. 
			Moreover, our method exhibits a small communication footprint and it achieves a comparably high solution accuracy with a limited number of iterations, whereby the local subproblems are of low computational complexity.
			We illustrate the performance of the proposed method on a problem from energy systems, i.e., we consider an optimal power flow problem with 708 buses. 
		\end{abstract}
		
		
	\end{frontmatter}
	
	\section{Introduction}
	Distributed and decentralized optimization algorithms are key for the optimal operation of networked systems.\footnote{\label{fn:dec} We refer to an optimization algorithm as being \textit{distributed} if one has to solve a (preferably cheap) coordination problem in a central entity/coordinator. We denote an optimization algorithm as being \textit{decentralized} in absence of such a coordinator and when the agents rely purely on neighbor-to-neighbor communication \citep{Bertsekas1989,Nedic2018}. 
		We call an algorithm \emph{essentially decentralized} if it has no central coordination but requires a small amount of central communication.
		Note that the definition of distributed and decentralized control differs~\citep{Scattolini2009}. }  
	Application examples range from  power systems \citep{Worthmann2015,Erseghe2015}, via optimal operation of gas networks networks \citep{Arnold2009}, to distributed control of data networks \citep{Low1999}.

	Classical distributed optimization algorithms used in the above works are, however, typically guaranteed to converge only for problems with convex constraints.
	Moreover, sufficiently accurate models are often non-linear leading to problems with non-convex constraints. 
	Thus, researchers either apply classical methods without convergence guarantees in a heuristic fashion \citep{Erseghe2015}, or they rely on simplified convex models \citep{Worthmann2015}.
	From an operations point of view, both approaches come with the risk of an unstable system operation.
	Moreover, classic approaches achieve linear convergence rates in the best case \citep{Hong2017,Yang2019}.

	\cite{Lu2018a}, \cite{Yan2011}, and \cite{Engelmann2019} propose distributed second-order methods with fast---i.e. superlinear---convergence guarantees for non-convex problems.
	These approaches rely on the exchange of quadratic models of the subproblems, which in turn implies a substantial need for communication and/or central coordination.
	In \citep{Engelmann2020c} we have shown how to overcome quadratic model exchange by a combination of active set methods and techniques from inexact Newton methods.
	However, in practice the detection of the correct active set is difficult and can  be  numerically unstable.

	In two recent papers, we have shown how to decompose interior point methods in an essentially decentralized fashion, i.e., decomposition is achieved without relying on any central computation \citep{Engelmann2021c,Engelmann2021}. 
	We do so by combining interior point methods with decentralized inner algorithms for solving the Newton step. 
	Interior point methods  have the advantage that an active set detection is avoided while  fast---i.e. superlinear---local convergence can be guaranteed for non-convex problems.	
	This note   considers the application of the essentially decentralized interior point method~(d-IP) to the optimal power flow problem which arises frequently  in power systems.

	\section{Problem Formulation}
	
	A common formulation of optimization problems in the context of networked systems is
	\begin{subequations} \label{eq:sepForm}
		\begin{align} 
			\min_{x_i,\dots,x_{|\mathcal{S}|}} \; \sum_{i \in \mathcal{S}} &f_i(x_i) \\
			\text{subject to }\quad  g_i(x_i)&=0, & \forall i \in \mathcal{S}, \label{eq:sepProbGi} \\
			h_i(x_i) &\leq 0, & \forall i \in \mathcal{S},\label{eq:sepProbHi} \\
			\sum_{i \in \mathcal{S}} A_ix_i &= b,\label{eq:consConstr}
		\end{align}
	\end{subequations}
	where, $\mathcal{S}=\{1,\dots, |\mathcal{S}|\}$ denotes a set of subsystems, each of which is equipped with an objective function $f_i:\mathbb{R}^{n_i} \rightarrow \mathbb{R}$ and equality and inequality constraints $g_i, h_i:\mathbb{R}^{n_i} \rightarrow \mathbb{R}^{n_{gi}},\mathbb{R}^{n_{hi}}$. 
	The matrices $A_i \in \mathbb{R}^{n_c \times n_i}$ and the vector $b\in \mathbb{R}^{n_c}$ are coupling constraints between the subsystems.
	
	\section{A Distributed Interior Point Method}
	Interior point methods reformulate problem \eqref{eq:sepForm} via a logarithmic barrier function and slack variables  $v_i\in \mathbb{R}^{n_{hi}}$, 	
		\begin{subequations} \label{eq:slackReform}
		\begin{align} 
			\min_{x_1,\dots,x_{|\mathcal{S}|},v_1,\dots,v_{|\mathcal{S}|}} \; \sum_{i \in \mathcal{S}} &f_i(x_i) - \mathds 1^\top \delta  \ln (v_i) \hspace{-1.5cm} \\
			\text{subject to }\quad  g_i(x_i)&=0, &  \forall i \in \mathcal{S}, \label{eq:eqCnstr}\\
			h_i(x_i) +v_i&= 0, \;\; v_i \geq 0,& \forall i \in \mathcal{S}, \label{eq:ineqCnstr}\\
			\sum_{i \in \mathcal{S}} A_ix_i &= b \label{eq:conCnstr}.
		\end{align}
	\end{subequations}
	The variable $\delta \in \mathbb{R}_+$ is a barrier parameter, $\mathds{1}= (1,\dots,1)^\top \in \mathbb{R}^{n_{hi}}$ and  the function $\ln(\cdot)$ is evaluated component-wise.
	Note that the inequality constraints are replaced by  barrier functions.
	Moreover, \eqref{eq:slackReform} and \eqref{eq:sepForm} share the same minimizers for $\delta \rightarrow 0$.
	
	The main idea of interior point methods is to solve \eqref{eq:slackReform} for a decreasing sequence of $\delta$.
	It is often too expensive to solve \eqref{eq:slackReform} to full accuracy---hence one typically performs a hand full Newton steps only \citep{Nocedal2006}.
	In this note we use a variant which computes only \emph{one} Newton step per iteration.

Next, we give a brief summary of distributed interior point methods; 
details are given in \citep{Engelmann2021c}.	

	\subsection{Decomposing the Newton Step}
An exact Newton step  $\nabla F^\delta(p)\Delta p = - F^\delta (p) $ applied to the first-order optimality conditions $F^\delta (p)=0$ of  \eqref{eq:slackReform} reads 
	\begin{align} \label{eq:structKKT}
		\begin{pmatrix}
			\hspace{-.3mm}\nabla F_1^\delta  &0& \dots& \tilde A_1^\top\hspace{-.3mm} \\
			0& \nabla F_2^\delta& \dots& \tilde A_2^\top  \hspace{-.3mm}\\
			\vdots & \vdots & \ddots& \vdots \\
			\tilde A_1 & \tilde A_2 & \dots&  0
		\end{pmatrix}
		\hspace{-1.9mm}
		\begin{pmatrix}
			\hspace{-.3mm}\Delta p_1 \hspace{-.3mm} \\
			\hspace{-.3mm}\Delta p_2\hspace{-.3mm} \\
			\vdots \\
			\Delta \lambda
		\end{pmatrix}
		\hspace{-1.4mm}
		=
		\hspace{-1.2mm} 
		\begin{pmatrix}
			-F_1^\delta \\
			-F_2^\delta   \\
			\vdots\\
			\hspace{-.6mm} b \hspace{-.3mm} -\hspace{-.5mm}  \sum_{i \in \mathcal{S}} \hspace{-.7mm} A_i x_i \hspace{-.2mm}
		\end{pmatrix},
	\end{align}
	where 
	\begin{align*}
		\nabla F_i^\delta =
		\begin{pmatrix}
			\nabla_{xx} L_i  & 0 & \nabla g_i(x_i)^\top  & \nabla h_i(x_i)^\top   \\
			0 &-V_i^{-1}  M_i & 0 & I  \\
			\nabla g_i(x_i) & 0 &0 & 0  \\
			\nabla h_i(x_i) & I & 0 & 0  \\
		\end{pmatrix} ,
	\end{align*}
	$M_i = \operatorname{diag}(\mu_i)$,
	and
	$
	\tilde A_i = 
	\begin{pmatrix}
		A_i &\; 0 &\; 0 & \;0
	\end{pmatrix}$,
cf. \citep[Thm. 12.1]{Nocedal2006}. 
Here,	 $p=(p_1,\dots,p_{|\mathcal S|},\lambda)$ and $p_i  = ( x_i, v_i, \gamma_i, \mu_i )$, where $\gamma_i$, $\mu_i$, and $\lambda$ are Lagrange multipliers assigned to \eqref{eq:eqCnstr}, \eqref{eq:ineqCnstr}, and \eqref{eq:conCnstr} respectively.
Note that the optimality conditions  $F^\delta (p)=0$ are parameterized by the barrier parameter $\delta$.
	
The coefficient matrix in \eqref{eq:structKKT} has an arrowhead structure which we exploit for decomposition. 
Note that each $\nabla F_i^\delta$ can be computed based on local information only. 
Assume  that  $\nabla F_i^\delta$ is invertible.
Then, one can reduce the KKT system \eqref{eq:structKKT} by  solving the first $S$ block-rows for $\Delta p_i$. 
Hence,  
\begin{align} \label{eq:delP}
	\Delta p_i = -  \left(\nabla F_i^\delta \right )^{-1}\left (F_i^\delta + \tilde A_i^\top \Delta \lambda  \right ) \text{ for all } i\in \mathcal{S}.
\end{align}
Inserting \eqref{eq:delP} into the last row of \eqref{eq:structKKT} yields 
\begin{equation} \label{eq:SchurComp}
	\begin{aligned}
		\Bigg (\sum_{i \in \mathcal{S}}\tilde A_i& \left(\nabla F_i^\delta \right )^{-1}\tilde A_i^\top \Bigg ) \Delta \lambda  \\
		&= \left(\sum_{i \in \mathcal{S}} A_i x_i -  \tilde A_i \left(\nabla F_i^\delta \right )^{-1}F_i^\delta  \right) - b. 
	\end{aligned}
\end{equation}
Define 
\begin{subequations} \label{eq:Schur}
	\begin{align} 
		S_i \doteq& \tilde A_i \left(\nabla F_i^\delta \right )^{-1}\tilde A_i^\top, \quad \text{and} \\
		s_i \doteq&  A_i x_i - \tilde A_i \left(\nabla F_i^\delta \right )^{-1}F_i^\delta -\dfrac{1}{|\mathcal{S}|} b.
	\end{align}
\end{subequations}
Then, equation \eqref{eq:SchurComp} is equivalent to
\begin{align} \label{eq:schurComp}
	\left (\sum_{i \in \mathcal{S}} S_i \right )\,\Delta \lambda - \sum_{i \in \mathcal{S}} s_i = S \Delta \lambda -s =0.
\end{align}

Observe that once \eqref{eq:schurComp} is solved, one can compute $\Delta p_1,\dots,\Delta p_{|\mathcal S|}$ locally in each subsystem based on $\Delta \lambda$ via back-substitution into \eqref{eq:delP}.
This way, we are able to solve~\eqref{eq:structKKT} in a hierarchically distributed fashion, i.e., we first compute $(S_i,s_i)$ locally for each subsystem and then collect $(S_i,s_i)$ in a coordinator.
One continues by solving~\eqref{eq:schurComp}  and  distributing $\Delta \lambda$ back to all subsystems $i \in \mathcal{S}$, which in turn use \eqref{eq:delP} to recover $\Delta p_i$.

\subsection{Decentralization}
Solving \eqref{eq:SchurComp} in a central coordinator is typically not preferred due to the large amount of information exchange for large-scale systems and due to safety reasons. 
Hence, we solve \eqref{eq:schurComp} in a decentralized fashion without central computation via decentralized inner algorithms. 

One can show that $S$ is symmetric and positive-semidefinite. 
Hence, one can apply a decentralized version of the conjugate gradient method~(d-CG) \citep{Engelmann2021}.
Alternatively, one can also use decentralized optimization algorithm by reformulating \eqref{eq:schurComp} as a convex optimization problem.

Typically it is expensive in terms of communication and computation to solve \eqref{eq:schurComp} to full accuracy by inner algorithms. 
Thus, we use techniques from inexact Newton methods to terminate the inner algorithms early based on the violation of the optimality conditions $F^\delta (p)=0$, cf. \cite[Chap. 7.1]{Nocedal2006} .
Doing so, one can save a large amount of inner iterations---especially in early outer iterations. 
When $\|F^\delta(p^k)\| $ gets closer to zero, we also force the residual of \eqref{eq:schurComp} to become smaller to guarantee convergence to a minimizer.

\subsubsection*{Updating Stepsize and the Barrier Parameter}
The barrier parameter $\delta$ and the  stepsize $\alpha$ in the for the Newton step $p^{k+1}= p^k + \alpha \Delta p^k$ requires a small amount of central communication but no central computation.
Indeed, it is possible to compute local surrogates $\{\alpha_i\}_{i\in \mathcal S}$ and $\{\delta_i\}_{i\in \mathcal S}$ and take their minimal/maximal values over all subsystems to obtain $(\alpha,\delta)$.

\subsubsection{The Overall Algorithm}
The overall distributed interior point algorithm is summarized in Algorithm~\ref{alg:d-IP}.
Algorithm~\ref{alg:d-IP} has local superlinear convergence guarantees for non-convex problems in case the barrier parameter and the residual in \eqref{eq:SchurComp}   decrease fast enough, cf. \cite[Thm. 2]{Engelmann2021c}.

\begin{algorithm}[t]
	\caption{Distributed Interior Point Method for  \eqref{eq:sepForm}}
	\begin{algorithmic}[1]
		\State Initialization: $p_i^0$ for all $i \in \mathcal{S}$, $\delta^0,\lambda^0$, $\epsilon$ \label{stp:1}
		\While{$\|F^0(p^k)\|_\infty > \epsilon$} 
		\State compute $(S_i^k,s_i^k)$ locally via \eqref{eq:Schur}  \label{stp:3}
		\While{residual of \eqref{eq:SchurComp} too large} \label{stp:4}
		\State iterate \eqref{eq:schurComp} via a  decentralized  algorithm \label{stp:5}
		\EndWhile\label{euclidendwhile}
		\State compute stepsize $\alpha^k$  and update $p^{k+1} = p^k + \alpha \Delta p^k$\label{stp:9}
		\State update $\delta^{k+1} < \delta^k$, $k \rightarrow k+1$
		\EndWhile\label{euclidendwhile}
		\State \textbf{return} $p^\star$
	\end{algorithmic} \label{alg:d-IP}
\end{algorithm}

	\section{Application to Optimal Power Flow} \label{sec:opf}
Optimal Power Flow~(OPF) problems aim at finding  optimal generator set-points in power systems while meeting  grid constraints and technical limits \citep{Frank2016}. 

The basic AC OPF problem reads
\begin{subequations}\label{eq:OPF}
	\begin{align} 
		&\min_{s,v \in \mathbb{C}^{N}} \;f(s) \\
		\text{subject to} \quad s-&s^d  = \operatorname{diag} (v) Y v^*,  \label{eq:PFeq} \\
		\underline p \leq \operatorname{re}(s) &\leq \bar p, \quad
		\underline q \leq \operatorname{im}(s) \leq \bar q, \label{eq:pwBounds}\\
		\underline v \leq \operatorname{abs}(v) &\leq \bar v, \quad v^1 = v^s. \label{eq:vBounds}
	\end{align}
\end{subequations}
Here,  $v \in \mathbb{C}^N$ are complex voltages, and $s \in \mathbb{C}^N$ are complex power injections   at all buses $N$.
The operators $\operatorname{re(\cdot)}$ and $\operatorname{im(\cdot)}$ denote the real part and imaginary part of a complex number, and $(\cdot)^*$ denotes the complex conjugate.
The objective function $f:\mathbb{C}^N\rightarrow \mathbb{R}$ encodes the cost of power generation.
The grid physics are described via the power flow equations \eqref{eq:PFeq}, where $Y \in \mathbb{C}^{N\times N}$ is the complex bus-admittance matrix describing  grid topology and parameters. 
Moreover, $s^d \in \mathbb{C}^N$ is a fixed  power demand.
The constraints \eqref{eq:pwBounds} describe technical limits on the power injection by generators, and \eqref{eq:vBounds} models voltage limits.
The second equation in \eqref{eq:vBounds} is a reference condition on the voltage at the first bus, $v^1$, where the complex voltage is constrained to a reference value~$v^s$.

Note that one can reformulate the OPF problem \eqref{eq:OPF} in form of \eqref{eq:sepForm} by introducing auxiliary variables. Different variants of doing do exist; here we rely on a reformulation from \cite{Muhlpfordt2020}.
	\subsection{A case study}

	\begin{figure}
		\centering
		\includegraphics[width=.85\linewidth,trim={6cm 11cm 5cm 11cm},clip]{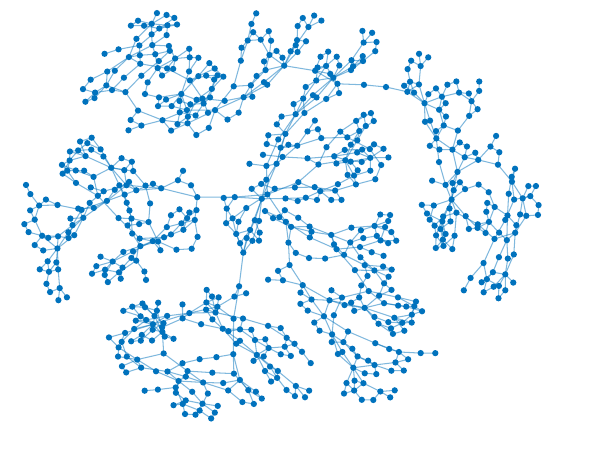}
		\caption{Six interconnected 118-bus systems.}
		\label{fig:grid}
	\end{figure}

	As a case study, we consider 6 interconnected IEEE 118-bus test systems shown in Fig.~\ref{fig:grid}.
	Each of these systems corresponds to one subsystem $i\in \mathcal{S}$ in problem~\eqref{eq:sepForm}.
	We use grid parameters from \texttt{MATPOWER}, and we interconnect the subsystems in an asymmetric fashion to generate non-zero flows at the interconnection points. 
	In total, we get an optimization problem with about $3.500$ decision variables. 

\begin{figure}
	\centering
	\includegraphics[width=1.08\linewidth]{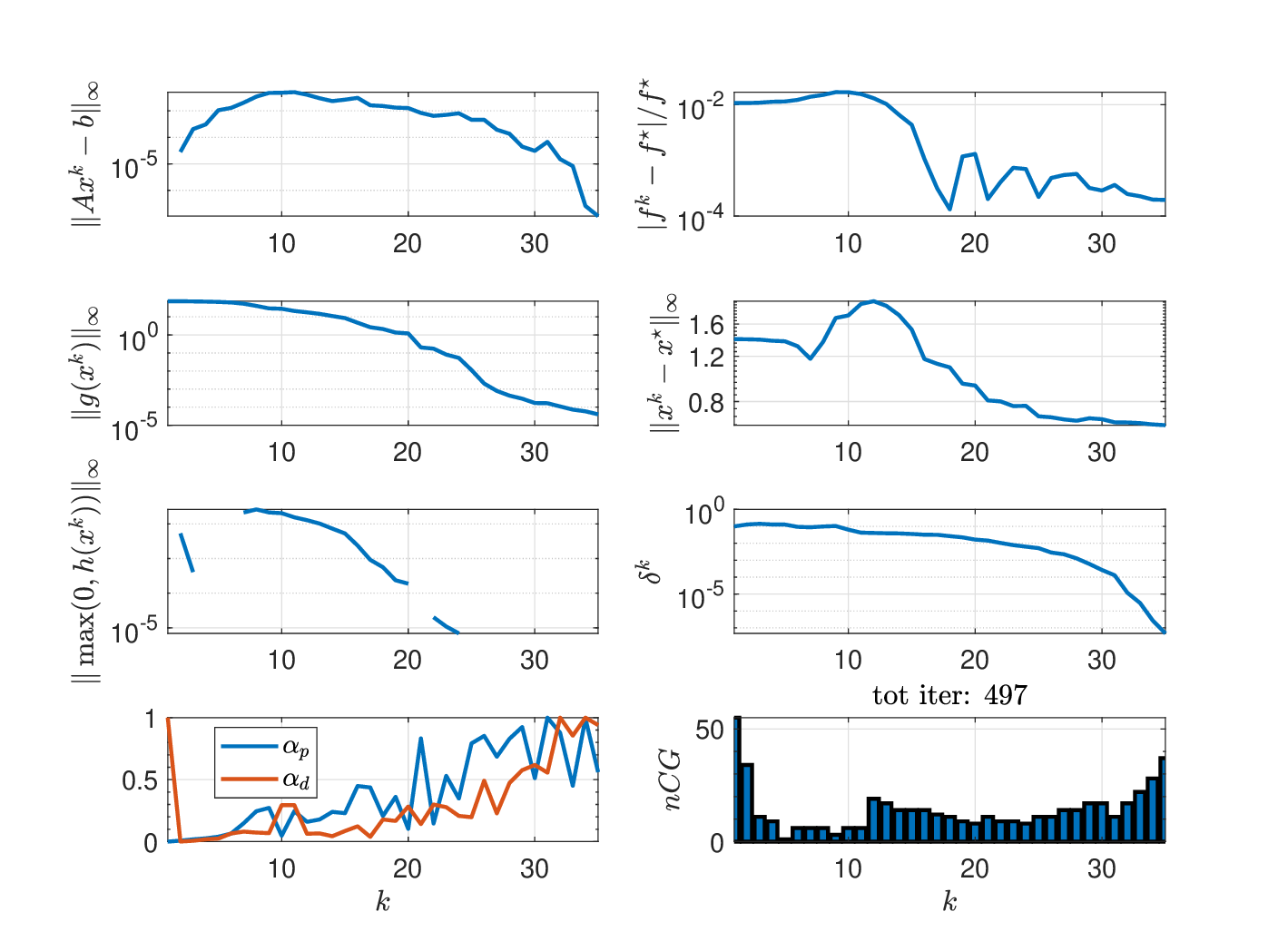}
	\caption{Convergence of Algorithm~\ref{alg:d-IP}.}
	\label{fig:dipconv}
\end{figure}

Fig.~\ref{fig:dipconv} depicts the convergence of Algorithm~\ref{alg:d-IP} over the iteration index $k$
with algorithm parameters from \citep{Engelmann2021c}.
The figure depicts the consensus violation $\|Ax^k-b\|_\infty$, which can be interpreted as the maximum mismatch of physical values at boundaries between subsystems.
Furthermore, the relative error in the objective  function $|f^k-f^\star|/f^\star$, the infeasibilities $\|g(x^k)\|_\infty$ and $\|\max(0,h(x^k))\|_\infty$, the distance to the minimizer $\|x^k-x^\star\|_\infty$, the number of inner iterations of d-CG, the barrier parameter sequence $\{\delta^k\}$, and the primal and dual step size $(\alpha^p,\alpha^d)$ are shown. 
The centralized solution $x^\star$ is computed via the open-source solver \texttt{IPOPT} \citep{Wachter2005}.

One  can observe that the consensus violation is at the level of $10^{-5}$ for all iterations.
This means that the iterates are feasible with respect to the power transmitted over transmission lines. 
This results from the fact that the consensus constraint \eqref{eq:conCnstr} is implicitly enforced when solving  \eqref{eq:schurComp} via d-CG.
A low consensus violation has the advantage, that one can terminate d-IP early and apply one local NLP iteration to obtain a feasible but possibly suboptimal solution.\footnote{Assuming that the local OPF problem is feasible for the current boundary value iterate.}
We note that feasibility is typically of much higher importance  than optimality in power systems, since feasibility ensures a safe system operation, cf. Remark~\ref{rem:feas}.
From  $\|g(x^k)\|_\infty$ and $\|\max(0,h(x^k))\|_\infty$\footnote{The blank spots in the plot for $\|\max(0,h(x^k))\|_\infty$ correspond to zero values, since $\log(0)=-\infty$.} in Fig.~\ref{fig:dipconv} one can see that feasibility is ensured to a high degree after 20-30 dIP iterations.
At the same time we reach a suboptimality level of almost $0.01 \%$, which is much smaller than in other works on distributed optimization for OPF, cf. \citep{Erseghe2015,Guo2017}.
Moreover, one can see that the distance to the minimizer $\|x^k-x^\star\|_\infty$ is still quite large due to the small sensitivity of $f$ with respect to the reactive power inputs. This is also common in OPF problems.

Regarding Algorithm~\ref{alg:d-IP} itself, one can see that the barrier parameter $\delta$ steadily decreases in each iteration. 
Moreover, during the first 20 iterations, comparably small step-sizes are used. The domain of local convergence is reached after around 30 iterations. Note that we use different stepsizes $\alpha_p$ for the primal variables and $\alpha_d$ for the dual variables.
Observe that due to the dynamic termination of inner d-CG iterations based on the inexact Newton theory, Algorithm~\ref{alg:d-IP} requires a small amount of inner iterations in the beginning and the number of iterations have to increase when coming closer to a local minimizer.
This saves a substantial amount of inner iterations compared to a fixed inner termination criterion.

The widely used Alternating Direction Method of Multipliers (ADMM) does not converge for the considered case. 
This seems to occur rarely, but was also reported in other works \citep{Christakou2017}.
Algorithm~\ref{alg:d-IP} requires 25 seconds for performing 35 iterations with serial execution.
The \texttt{MATPOWER} solver \texttt{MIPS} needs about 13 seconds when applied to the distributed formulation and 2 seconds when applied to the centralized problem formulation.
Executing 497 ADMM iterations---this reflects the number of d-CG iterations in Algorithm~\ref{alg:d-IP}--- requires 210 seconds with serial execution.
This illustrates the large computation overhead of ADMM in the local steps, since here one has to solve an NLP in each iteration and for each subsystem. 
In contrast, d-IP only needs to perform one matrix inversion every outer iteration.
All simulations are performed on a standard state-of-the-art notebook.

\begin{remark}[{Sufficient feasibility in power systems }] \label{rem:feas}
	
Note that feasibility in a range of $10^{-3}$ to $10^{-5}$ is typically sufficient for a safe power system operation.
The  parameters in the OPF problem \eqref{eq:OPF}, such as power demands and line parameters,  induce uncertainty to the problem, which is typically much larger than this level \citep{Kim2000}. 
Hence, in applications there is typically little-to-no benefit in solving OPF problems to machine precision.
\end{remark}

	\section{Summary \& Outlook}
	We have presented an essentially decentralized interior point method for distributed optimization in energy networks with advantageous properties in terms of convergence guarantees, communication footprint, and practical convergence.
	We have illustrated the performance of our method on a 708-bus case study.
Future work will consider  improvements in  implementation aspects of d-IP, where we aim faster execution times and at scalability up to several thousand buses. 
	
	\renewcommand*{\bibfont}{\footnotesize}
	\bibliography{./paper}               
	
	
	
	
	
	
	
	
\end{document}